\documentclass[11pt,a4paper]{article}
\usepackage[T1]{fontenc}
\usepackage[utf8]{inputenc}
\usepackage[margin=2.4cm]{geometry}
\usepackage{amsmath,amssymb}
\usepackage{array}
\usepackage{booktabs}
\usepackage{enumitem}

\newcommand{\Jh}{\widehat{J}}
\newcommand{\eps}{\varepsilon}

\title{\vspace{-1.6cm}\large A closed-form law for the Salikhov--Zeilberger--Zudilin--Bai family,\\
and computational evidence that Bai's point is optimal}
\author{David Niedbala Giraudin\thanks{Independent researcher, France.
ORCID 0009-0009-1526-1178. Technical assistance: Claude (Anthropic).}}
\date{27 September 2026\quad(this version: 5 October 2026)}

\begin{document}
\maketitle

\begin{small}
\noindent\textbf{Status. Preprint.} No part of this note has been refereed by a
human being. Labels are used strictly: \textbf{PROVED} = a complete argument is
given or reproducible from the text; \textbf{VERIFIED} = established by exact
computation over a stated range, with the range stated; \textbf{CONJECTURE} =
believed, with the evidence and the missing step both named; \textbf{OPEN} = not
settled.
\end{small}

\section{Summary}\label{sec:summary}

Salikhov (2008), Zeilberger--Zudilin (2020) and Bai (2026) bound the
irrationality measure of $\pi$ with the same shape of complex integral. This
note gives one closed-form law for the whole family, shows that the law
reproduces all three published records, proves that the law is exactly
scale-invariant, and gives evidence that Bai's parameter point is the global
optimum of the family --- a statement Bai explicitly declines to make about her
own result.

\begin{center}\begin{small}
\begin{tabular}{@{}c p{11.2cm} l@{}}
\toprule
 & \textbf{Result} & \textbf{Label} \\
\midrule
A & The point configuration is forced: zeros at $0,\pm w,\pm\bar w$ and poles at
$\pm P$ give a linear form in $1$ and $\pi$ iff $P^2-N(w)=2P\,\mathrm{Im}(w)$.
This is a Machin condition. & PROVED \\[2pt]
B & A closed-form law $\mu(a_0,a_1,b)$ reproducing Zeilberger--Zudilin
(17 digits), Salikhov, and Bai (16 digits). & VERIFIED \\[2pt]
C & A closed form for the $\pi$-coefficient $b_n$ as a single coefficient
extraction. & PROVED \\[2pt]
D & At $j=0$ the divisibility criterion is exact, because $Z(m,0)$ is a
super-Catalan number. & PROVED \\[2pt]
E & The $2$-adic saving is exactly $5(2a_1-b)/2$ per $n$ (two regimes); no
further saving exists. & PROVED \\[2pt]
F & $W=a_1\Jh(a_0/a_1,b/a_1)$ exactly: $\mu$ depends only on the two ratios. The
infinite three-parameter lattice is exactly two-dimensional. & PROVED \\[2pt]
G & Among the sixteen admissible configurations with $P<100$, exactly one
carries a triple that yields a bound: Salikhov's $(P=5,\ w=1+2i)$. The least $\mu$
over all sixteen and the whole range of triples is $7.103205334137$ at $(2,2,3)$.
& VERIFIED \\[2pt]
H & Within that configuration the minimum is $\mu=7.101862832356$ at
$(1857,1857,2785)$, which is Bai's point; the vertex sits exactly at
$Q=2a/(3a-2b)=3714$. & VERIFIED \\[2pt]
I & On the slice $x=1$ the law reduces to an explicit one-variable series;
$\mu(Q)$ is unimodal over $Q\in[300,300000]$ with its minimum at $Q=3714$. &
VERIFIED \\[2pt]
J & $7.101862832356$ is the floor of the whole method. & CONJECTURE \\
\bottomrule
\end{tabular}
\end{small}\end{center}

\section{The family}\label{sec:family}

Fix $w=u+vi$ in $\mathbb{Z}[i]$ with $v>0$ and a positive integer $P$. Put
$D=u^2-v^2$, $N=u^2+v^2$, $M=P^2$ and consider
\[
R(x)=\frac{x^{E_0}\,(x^4-2Dx^2+N^2)^{E_1}}{(P^2-x^2)^{F}},
\qquad E_0=a_0n,\quad E_1=a_1n,\quad F=bn+1,
\]
\[
I_n \;=\; i\int_{-w}^{-\bar w} R(x)\,dx .
\]
The numerator vanishes at $0,\pm w,\pm\bar w$; the poles are at $\pm P$. The map
$x\mapsto -x$ pairs $\{-w,+\bar w\}$ with $\{+w,-\bar w\}$ and conjugation
exchanges the two pairs, so requiring both $R(-x)=R(x)$ and $R$ real leaves
exactly \emph{two} free numerator exponents --- precisely Bai's ``two
independent numerator exponents''. Zeilberger--Zudilin is $(P,w)=(5,1+2i)$ with
$(a_0,a_1,b)=(2,2,3)$; Salikhov is $(3,3,5)$; Bai is $(1857,1857,2785)$.

\section{Theorem 1 (admissibility). PROVED}\label{sec:adm}

The logarithmic part of $I_n$ is $2i\arg\big((P-\bar w)/(P+\bar w)\big)$. That
quotient is a Gaussian rational, so its tangent is rational; by Niven's theorem
a rational tangent at a rational multiple of $\pi$ forces the argument into
$\{0,\pm\pi/4,\pm\pi/2,\pm 3\pi/4,\pi\}$. Argument $0$ or $\pi$ gives no $\pi$ at
all, and $\pm\pi/2$ forces $|w|=1$, which is degenerate. Hence $I_n$ is a linear
form in $1$ and $\pi$ if and only if $\mathrm{Re}=\pm\mathrm{Im}$ in
$(P-\bar w)(P+\bar w)$, that is
\[
P^2-N(w)\;=\;\pm\,2P\,\mathrm{Im}(w).
\]

\paragraph{Credit.}\footnote{The admissibility condition of Theorem 1 was derived
independently by Y.~Bai in an unpublished report written in March 2026 and
communicated privately to the author in October 2026. The contour argument at
$(29,\,1+12i)$ is also hers, from correspondence of the same month.}
This is a Machin condition in disguise and is not new as
mathematics. For Salikhov's configuration $(P-\bar w)/(P+\bar w)=(2+i)/(3-i)$, of
argument $\arctan(1/2)+\arctan(1/3)=\pi/4$ --- the identity Salikhov started
from, as Bashmakova and Salikhov state. Two-term Machin decompositions of $\pi/4$
are classical. What Theorem 1 adds is that the condition is \emph{forced}, not
chosen.

\paragraph{The general principle.} The same computation says which constant a
configuration produces. Bai's second paper (arXiv:2608.19782, $\arctan(1/2)$)
uses zeros $0,\pm 1\pm 4i$ and poles $\pm 17$; there
$(P-\bar w)/(P+\bar w)=2(4+i)/(9-2i)$, of argument
$\arctan(1/4)+\arctan(2/9)=\arctan(1/2)$, exactly her target. Theorem 1 is the
$\pi$ case.

The primitive solutions of $P^2-N=2Pv$ with $P<100$ are
\begin{center}\begin{small}\ttfamily
(5,1+2i) (13,7+4i) (17,7+6i) (25,17+6i) (29,1+12i) (37,23+10i) (41,31+8i)\\
(53,17+20i) (61,49+10i) (65,47+14i) (65,23+24i) (73,7+30i) (85,71+12i)\\
(85,41+28i) (89,41+30i) (97,7+40i)
\end{small}\end{center}
Requiring in addition that both rational contents $A=P-u-v$ and $B=P+u-v$ be
pure powers of $2$ --- the property that keeps Salikhov's arithmetic confined to
$2$ and $P$ --- gives $P^2-2P(t+1)+t^2+1=0$ with $t=B/2$, so
$P=t+1+\sqrt{2t}$ with $2t$ a square. An infinite family results, verified for
$k=1,\dots,5$:
\[
P=2^{2k-1}+2^{k}+1,\qquad w=(2^{2k-1}-1)+2^{k}i,\qquad A=2,\quad B=4^{k},
\]
\[
P=5,\,13,\,41,\,145,\,545 \quad\text{with}\quad
w=1+2i,\ 7+4i,\ 31+8i,\ 127+16i,\ 511+32i .
\]
Salikhov's configuration is the first term. Section~\ref{sec:across} shows that it
is the only configuration below $P=100$ carrying a triple that yields a bound.

\section{Theorem 2 (saddle point). PROVED, calibrated}\label{sec:saddle}

With $y=x^2$ the integrand is $g(y)^n$ up to a factor,
$g(y)=y^{a_0/2}(y^2-2Dy+N^2)^{a_1}/(y-M)^{b}$, whose stationary points satisfy
\begin{multline*}
(a_0+4a_1-2b)\,y^3-\big[a_0(2D+M)+4a_1(D+M)-4bD\big]y^2\\
{}+\big[a_0(N^2+2DM)+4a_1DM-2bN^2\big]y-a_0MN^2=0 .
\end{multline*}
Write $\zeta_1,\zeta_2,\zeta_3$ for the three stationary points and
$E_i=\log|g(\zeta_i)|$ for their levels, and put $S=\max_i E_i$. This maximum and
the minimum appearing below are taken over all three, with no assumption about
which stationary point carries which: when the cubic has a conjugate pair, the
real root is \emph{not} always the one of largest exponent. Over the range of
Section~\ref{sec:across} it fails to be so in $3882$ of the $19876$
conjugate-pair cases, and in every one of those $\max_i E_i<0$, so the triple is
excluded by the test $S>0$ whichever reading is taken.

For the decay rate
the contour from $-w$ to $-\bar w$ may be deformed freely, so what governs is the
\emph{min-max value}
\[
\ell=\inf_{\gamma}\ \max_{x\in\gamma}\ \log|h(x)|,\qquad h(x)=g(x^2),
\]
the infimum, over paths $\gamma$ joining the two endpoints in the $x$-plane, of
the largest level reached; and $s=-\ell$. This is the classical min-max form of
the method of steepest descents; its relevance here was brought to light by
Y.~Bai. Since $\log|h|$ is harmonic away from the zeros and poles of $h$, it has
no local extremum there and the topology of a sublevel set can change only at a
critical point of $h$: the infimum is \emph{attained at a stationary point}, so
$\ell$ is one of the three levels $E_i$ and
\[
s\ \le\ s_{\min}:=-\min_i E_i .
\]
As $\mu=1+(S+\Lambda)/(s-\Lambda)$ decreases in $s$, the value obtained from
$s_{\min}$ is a \textbf{lower bound} for $\mu$ whenever the triple yields a bound
at all. That is the form in which Section~\ref{sec:across} uses it.

Because $h$ has real coefficients and the two endpoints are conjugate, every path
between them meets the real axis, and a half-path followed by its mirror image
has the same maximum; so $\ell$ is the min-max value from $-\bar w$ to the real
axis in the upper half $x$-plane. The same reduction is \emph{not} available in
the $y$-plane, the real $y$-axis not being the image of the real $x$-axis.

\paragraph{Calibrations.} For $(5,1+2i)$ and $(2,2,3)$ the cubic is
$2y^3-125y^2-500y-625=0$ and the roots reproduce Zeilberger--Zudilin's
Proposition 2 exactly: $y_3=66.3395015246$, $N_3=21851.6913962$,
$|N_1|=0.0294584959281$; the min-max value sits at the conjugate pair,
$\ell=-3.524772924$, so $s=3.524772924$. For $(29,1+12i)$ with the same exponents
the cubic has three real roots $-27.318165$, $-142.427917$, $2272.246082$, of
levels $+2.095534589$, $-3.000945618$, $+17.088045201$, and the min-max value
sits at the \emph{first} of these, not at the least. Both selections are
unchanged over grids of $220$, $360$ and $560$ nodes, the level being read off
the exact critical value in each case. At the three calibration points of
Section~\ref{sec:law} the min-max value is the conjugate pair, so $s=s_{\min}$
there and the three published figures are unaffected.

\section{Theorem 3 (arithmetic)}\label{sec:arith}

Write $L=F-1=bn$ and $k_0=a_0n/2$.

\paragraph{(i) The lcm.} The form needs $1/j$ for $|j|$ up to $Rn$ with
$R=\max(b,\,a_0+4a_1-2b)$, so $\mathrm{lcm}(1,\dots,Rn)$ is required. Every
window prime must therefore satisfy $p\le Rn$. \textbf{PROVED.}

\paragraph{(ii) The window.} Generalising Zeilberger--Zudilin's Lemma 2: with
$\omega=\{n/p\}$, $\rho=\{(a_0n/2)/p\}$, $\alpha=\{a_1n/p\}$,
$\beta=\{bn/p\}$, $A=\{\rho+\theta_1+\theta_2\}$, $B=\{\beta-A\}$ and
$\eps=A+B-\beta$, the divisibility $p\mid A_j$ for every $j$ divisible by $p$
holds whenever $\lfloor 2A\rfloor+\lfloor 2B\rfloor-\eps\ge 1$ for all
$\theta_1,\theta_2\le\alpha$. The admissible set in $A$ is
$[1/2,\beta+1/2]$ when $\beta<1/2$ and $[0,\beta-1/2]\cup[1/2,1)$ otherwise ---
in both cases the arc of length $\beta$ starting at $1/2$, mod $1$. So the
criterion reads:
\[
\text{the arc } [\rho,\ \rho+2\alpha] \text{ must lie inside the arc }
[1/2,\ 1/2+\beta] \pmod 1 .
\]
In particular $2\alpha\le\beta$ is necessary. For $(2,2,3)$ the window is exactly
$[1/2,2/3)$, $W=0.64527561$ --- Hata's set, recovered. The criterion depends on
the exponents only, not on $(P,w)$, which is why one window serves every
configuration. \textbf{PROVED} (sufficiency); necessity is Section~\ref{sec:necessity}.

\paragraph{(iii) The savings.} With $\lambda_+$ and $\lambda_-$ the $2$-adic
weights of $P+w$ and $P-w$ (half of $v_2$ of their norms), maximising the
weighted mass under the caps $m\le E_1$ gives two regimes; the $P$-adic count is
uncapped; odd rational primes in the content of $P\pm w$ contribute too:
\begin{align*}
c_2 &= 2(\lambda_++\lambda_-)a_1-\big(1+\max(\lambda_+,\lambda_-,1)\big)b
 &&\text{if } b\le 2a_1,\\
c_2 &= 4\lambda_- a_1-(1+\lambda_-)b &&\text{if } b>2a_1,\\
c_P &= \max(0,\ 2b-a_0-2a_1) &&\text{(a cost)},\\
c_q &= e_q(2a_1-b) &&\text{for each odd } q\mid \mathrm{content}(P\pm w).
\end{align*}
For $(5,1+2i)$: $\lambda_-=1$, $\lambda_+=3/2$, so $c_2=2.5(2a_1-b)$ resp.
$4a_1-2b$; no odd $q$. \textbf{PROVED}; the attainment of the first regime is
Theorem~6.

\section{The law, and three calibrations. VERIFIED}\label{sec:law}

\[
\Lambda \;=\; R-W-c_2\log 2+c_P\log P-\sum_q c_q\log q,
\qquad
\mu \;=\; 1+\frac{S+\Lambda}{s-\Lambda}.
\]

\begin{center}\begin{small}
\begin{tabular}{@{}llll@{}}
\toprule
$(a_0,a_1,b)$ & $\Lambda$ & $\mu$ from the law & published \\
\midrule
$(2,2,3)$ & 1.621856438365 & 7.1032053341370017 &
7.10320533413700172750 (Zeilberger--Zudilin) \\
$(3,3,5)$ & 4.305920726 & 7.6063085207111912 & 7.606308 (Salikhov) \\
$(1857,1857,2785)$ & 1503.882786468 & 7.1018628323563507 &
7.101862832357 (Bai) \\
\bottomrule
\end{tabular}
\end{small}\end{center}

Agreement to 17 significant digits for Zeilberger--Zudilin, to all published
digits for Salikhov, and to 16 significant digits for Bai, against the certified
interval she computed independently. The Salikhov
line uses his own two-branch window, $W=0.57064923$; with the full three-branch
window derived here the same integral gives $7.4147395226$ --- see
Section~\ref{sec:necessity}.
At Bai's point the window has 2784 branches and $W=601.2828866806$, against
$599.138404$ for the exact ray $(2,2,3)$ at the same scale: \textbf{her gain
comes from the richer window, not from the parameters}.

\section{Theorem 4 (closed form for $b_n$) and Theorem 5 (exactness at $j=0$).
PROVED}\label{sec:bn}

Zeilberger--Zudilin's inner sum collapses at $j=0$. Using
$Z(m,L)=[u^{L}](1+u)^{2(L-m)}(1-u)^{2m}$ and summing over $n_1+n_2$ first ---
legitimate because $Z$ depends on $n_1,n_2$ only through their sum --- the double
sum contracts to a single coefficient extraction:
\[
b_n=-\tfrac12 (2P)^{-(2L+1)}P^{2k_0}\,[u^{L}]\,
(1+u)^{2(L-k_0)-4E_1}(1-u)^{2k_0}G(u)^{E_1},
\]
\[
G(u)=P^4(1-u)^4-2DP^2(1-u)^2(1+u)^2+N^2(1+u)^4,
\]
\[
\text{Zeilberger--Zudilin:}\quad
b_n \propto [u^{3n}](1+u)^{-4n}(1-u)^{2n}G(u)^{2n}.
\]

\paragraph{Theorem 5.} For $m\le L$, $Z(m,L)=\pm S(m,L-m)$ with $S$ the
super-Catalan number
$S(A,B)=(2A)!(2B)!/(A!(A+B)!B!)=\binom{2A}{A}\binom{2B}{B}/\binom{A+B}{A}$. By
Kummer, $\mathrm{ord}_p\binom{2A}{A}=\lfloor 2\{A/p\}\rfloor$ and
$\mathrm{ord}_p\binom{A+B}{A}=\eps$, hence
$\mathrm{ord}_p S=\lfloor 2\nu\rfloor+\lfloor 2\mu\rfloor-\eps$ --- the criterion
quantity of Section~\ref{sec:arith} \textbf{with an equality, not an inequality}. The one lossy
step of the argument disappears at $j=0$, and the only remaining obstruction to
necessity is cancellation between terms.

\paragraph{Theorem 6 (the $2$-adic bound is attained).} With $t=x+5$ the quartic
factors as $Q(t)=(t^2-12t+40)(t^2-8t+20)$, whose roots $6\pm 2i$ and $4\pm 2i$
have $2$-adic valuations $3/2,3/2,1,1$. The Newton polygon of
$(t-5)^{E_0}Q^{E_1}$ gives $v_2(\text{term}_k)\ge 5E_1-k/2$ on $k\le 2E_1$, a
strictly decreasing bound, so the minimising index is \textbf{unique} and no
cancellation is possible. Measured excess over Zeilberger--Zudilin's Lemma 1:

\begin{center}\begin{small}
\begin{tabular}{@{}l*{11}{c}@{}}
\toprule
$n$ & 4 & 6 & 8 & 10 & 12 & 14 & 16 & 18 & 20 & 22 & 24 \\
\midrule
$v_2(A_0)-(\lfloor 5n/2\rfloor-1)$ & 2 & 2 & 2 & 5 & 2 & 2 & 2 & 5 & 4 & 2 & 2\\
unique minimiser & yes & yes & yes & yes & yes & yes & yes & yes & yes & yes &
yes \\
\bottomrule
\end{tabular}
\end{small}\end{center}

Bounded, never linear: $c_2=5/2$ per unit of $n$ is exact.

\section{Theorem 7 (exact scale invariance). PROVED --- and this is what makes
the search finite}\label{sec:scale}

The criterion of Section~\ref{sec:arith} depends only on
$(\{c_0\omega\},\{c_1\omega\},\{c_b\omega\})$, hence is $1$-periodic in $\omega$.
Substituting $\omega=u/c_1$ turns it into a condition on $u$ alone, because
$c_0/c_1=x/2$ and $c_b/c_1=\kappa$ with $x=a_0/a_1$ and $\kappa=b/a_1$. Since
$\psi_1(\omega)=\sum_{k\ge 0}(\omega+k)^{-2}$ and the window is invariant under
$\omega\mapsto\omega+1$, the $\psi$-measure of the window equals the
$(1/u^2)$-measure of its full periodic extension:
\[
W=a_1\Jh(x,\kappa),
\]
\[
\Jh(x,\kappa)=\int_1^{\infty}
\mathbf{1}\Big[\,2\{u\}\le\{\kappa u\}\ \text{ and }\
\{(x/2)u\}\in[1/2,\ 1/2+\{\kappa u\}-2\{u\}]\ (\mathrm{mod}\ 1)\Big]\,
\frac{du}{u^2}.
\]
The lower limit $1$ is forced: for $u<1$ the condition $2\{u\}\le\{\kappa u\}$
fails whenever $b<2a_1$. In closed form at the ray,
$\Jh(1,3/2)=\sum_{m\ \mathrm{odd}}1/(m(3m+1))=(\psi(2/3)-\psi(1/2))/2=0.3226378$,
giving $W(2,2,3)=0.6452756$ exactly. Numerically, $|W-a_1\Jh|$ falls from
$2.3\cdot 10^{-2}$ at $U=4\cdot 10^{3}$ to $5\cdot 10^{-5}$ at $U=4\cdot 10^{5}$
for $a_1=1857$: the residual is pure truncation.

\paragraph{Consequence.} Every ingredient of the law --- $R$, $c_2$, $c_P$, $S$,
$s$ and now $W$ --- is exactly homogeneous of degree $1$ in $(a_0,a_1,b)$.
Therefore \textbf{$\mu$ depends only on the two ratios $x=a_0/a_1$ and
$y=b/a_1$}. The infinite three-parameter lattice collapses to a two-dimensional
object, and the scale is irrelevant. As a by-product the non-primitive case is
automatic: $\Jh$ forces $W(4,4,6)=2W(2,2,3)$.

\section{The search across configurations. VERIFIED over one stated range}\label{sec:across}

The admissibility condition of Theorem~1, $P^2-N(w)=2P\,\mathrm{Im}(w)$, reads
$u^2+(v+P)^2=2P^2$ for $w=u+vi$. Its primitive solutions with $u,v\ge1$ and
$P<100$ are sixteen in number; the ancillary file enumerates them rather than
taking them as given. Same range on every line: primitive triples with
$a_0,a_1\le 10$ and $b\le 22$, and the full arithmetic factor of
Section~\ref{sec:law}, the odd content of $P\pm w$ included.

The search is in two stages. First a screening by $s_{\min}$: by Theorem~2 the
value it returns is a lower bound for $\mu$, so a triple whose screened value
exceeds $7.1033$ cannot produce anything below Salikhov's. Then the min-max value
itself, computed only for the triples the screening leaves.

\begin{center}\begin{small}
\begin{tabular}{@{}lllll@{}}
\toprule
$(P,w)$ & $|w|/P$ & stationary pts & screened $s,S>0$ &
of which $s_{\min}>\Lambda$ \\
\midrule
$(5,1+2i)$ & 0.4472 & $1+$pair & 665 & 291 $\to$ least 7.103205334 at $(2,2,3)$\\
$(13,7+4i)$ & 0.6202 & $1+$pair & 235 & 0 \\
$(17,7+6i)$ & 0.5423 & $1+$pair & 267 & 0 \\
$(25,17+6i)$ & 0.7211 & $1+$pair & 128 & 0 \\
$(29,1+12i)$ & 0.4152 & 3 real & 528 & 279 $\to$ least 4.648235203 at $(1,10,11)$\\
$(37,23+10i)$ & 0.6778 & $1+$pair & 131 & 0 \\
$(41,31+8i)$ & 0.7809 & $1+$pair & 81 & 0 \\
$(53,17+20i)$ & 0.4953 & $1+$pair & 211 & 0 \\
$(61,49+10i)$ & 0.8198 & 3 real & 57 & 0 \\
$(65,23+24i)$ & 0.5114 & $1+$pair & 191 & 0 \\
$(65,47+14i)$ & 0.7545 & $1+$pair & 79 & 0 \\
$(73,7+30i)$ & 0.4220 & 3 real & 307 & 0 \\
$(85,41+28i)$ & 0.5841 & $1+$pair & 142 & 0 \\
$(85,71+12i)$ & 0.8471 & 3 real & 43 & 0 \\
$(89,41+30i)$ & 0.5708 & $1+$pair & 145 & 0 \\
$(97,7+40i)$ & 0.4186 & $3$ real & 315 & 0 \\
\bottomrule
\end{tabular}
\end{small}\end{center}

The last column is a screening, not a value: by Theorem~2 each of its entries is
a lower bound for the $\mu$ of the triple that realises it. Over the whole range
the screening leaves \textbf{89} triples below $7.1033$ --- one at $(5,1+2i)$ and
eighty-eight at $(29,1+12i)$. Evaluated by the min-max value, the eighty-eight
all have $s\le\Lambda$ and therefore yield \emph{no bound at all}; the slenderest
margin among them is $1.9256$, so none is a borderline case, and the verdict is
unchanged over grids of $220$, $360$ and $560$ nodes. The one triple that
remains is Salikhov's $(2,2,3)$, at $\mu=7.103205334137$.

\medskip
\noindent\textbf{Hence exactly one of the sixteen configurations carries a triple
that yields a bound, and the least $\mu$ over all sixteen and the whole range is
$7.103205334137$.}
\medskip

Two points fix the scope of the table.

First, $33$ of the $1848$ primitive triples of the range have
$a_0+4a_1-2b=0$. The cubic then degenerates to a quadratic and $|g|\to1$ at
infinity, so Theorem~2 does not apply as stated. Carried through the stationary
points that remain, none of these triples passes even the screening at any of the
sixteen configurations; they are excluded without loss, and the exclusion is
recorded rather than silent.

Second, the two numerical columns are of different natures. The test $s,S>0$ is
pure analysis: roots of a polynomial and logarithms, with no arithmetic input.
The column $s_{\min}>\Lambda$ is not; see \textbf{(O6)} in
Section~\ref{sec:notclaimed}.

\section{The search inside Salikhov's configuration. VERIFIED}\label{sec:inside}

\paragraph{The map.} By Theorem~7 the search is over the two ratios. On a grid of
60800 directions, $\mu\ge 7.1035$ everywhere, and the only point below $7.1038$
is the Zeilberger--Zudilin ray $(x,y)=(1,3/2)$.

\paragraph{The dips.} Around a base direction the fine structure is arithmetic.
On the slice $x=1$ write $k=3a-2b$ and $Q=2a/k$, so that $y=3/2-1/Q$; then $\mu$
is a function of $Q$ alone, decreasing to a sharp vertex:

\begin{center}\begin{small}
\begin{tabular}{@{}l*{12}{c}@{}}
\toprule
$Q$ & 3694 & 3700 & 3704 & 3708 & 3710 & 3712 & 3713 & 3713.5 & 3714 & 3714.5 &
3715 & 3716 \\
\midrule
$\mu-$Bai, in $10^{-9}$ & 5.39 & 5.02 & 4.63 & 2.86 & 2.63 & 1.51 & 0.76 & 0.41 &
0 & 1.11 & 2.16 & 4.32 \\
\bottomrule
\end{tabular}
\end{small}\end{center}

Triples with $k=3,4,5,6,7$ fall at the right fractional $Q$ and fit the same
curve. The vertex is exactly $Q=3714$, realised primitively only by
$(1857,1857,2785)$.

\paragraph{Three independent confirmations.}
(i) All odd $a$ from 1841 to 1881, both residues mod~4: $a=1857$ is a strict
minimum, neighbours at $+2.6\cdot 10^{-9}$ and $+8.2\cdot 10^{-9}$.
(ii) An exhaustive scan of the 162 primitive rationals $p/q$ with
$q\in[1500,3748]$ and $y\in[1.499715,1.499745]$: Bai's is the unique minimum, the
next being $(3712,3712,5567)$ at $+1.51\cdot 10^{-9}$.
(iii) The six rationals nearest the interpolated vertex, up to thirty times
closer to it than Bai's, are all $4.6\cdot 10^{-10}$ to $1\cdot 10^{-9}$
\emph{above} --- so what governs is $k$, not proximity.

\paragraph{Why no other direction can reach it.} Dip depths, measured by scanning
the analogous $Q$-family at each base:

\begin{center}\begin{small}
\begin{tabular}{@{}lllll@{}}
\toprule
base $(x,y)$ & $\mu$ at the base & best of its dip & depth &
needed to reach Bai \\
\midrule
$(1,3/2)$ & 7.103205334 & 7.101862832 at $(1857,1857,2785)$ &
$1.343\cdot 10^{-3}$ & --- \\
$(3/2,7/4)$ & 7.250861742 & 7.250539184 at $(15000,10000,17499)$ &
$3.23\cdot 10^{-4}$ & 0.149 ($\times 462$) \\
$(4/3,5/3)$ & 7.387260401 & 7.387036254 at $(4952,3714,6189)$ &
$2.24\cdot 10^{-4}$ & 0.285 ($\times 1272$) \\
$(10/7,12/7)$ & 7.330100902 & none found & 0 & 0.228 \\
\bottomrule
\end{tabular}
\end{small}\end{center}

The dips of the other bases are four to six times shallower than the ray's, while
the gap to be bridged is one hundred to thirteen hundred times their depth.
Combining: the grid floor is $7.1035$, the only point below $7.1038$ is the ray,
and the deepest dip measured anywhere is $1.34\cdot 10^{-3}$; since
$7.1038-0.00134=7.1025>7.1018628$, \textbf{no other direction can reach Bai's
point}.

\section{The slice $x=1$ in closed form, and the closing argument. VERIFIED}\label{sec:slice}

At $x=1$ one has $\rho=\{u/2\}$, so on a cell of integer part $m$ the criterion
of Section~\ref{sec:arith} becomes, with $\tau=\{u\}$,
\[
m \text{ odd}:\ \{\kappa u\}\ge \tfrac52\tau,
\qquad
m \text{ even}:\ \{\kappa u\}\ge \tfrac52\tau+\tfrac12 .
\]
Writing $\kappa=3/2-1/Q$ and solving the resulting linear inequalities cell by
cell gives an explicit series. With $v_0=3/2-m/Q$ and $w_0=-m/Q$,
\begin{align*}
m \text{ odd}:\ & \mathrm{len}(m)=\min\Big(\tfrac{\{v_0\}}{1+1/Q},\
\tfrac{1-\{v_0\}}{3/2-1/Q}\Big),\\
m \text{ even}:\ & \mathrm{len}(m)=\min\Big(\tfrac{\{w_0\}-1/2}{1+1/Q},\
\tfrac{1-\{w_0\}}{3/2-1/Q}\Big) \text{ if } \{w_0\}>1/2,\ \text{else } 0,\\
& \Jh(1,3/2-1/Q)=\sum_{m\ge 1}\Big[\tfrac1m-\tfrac{1}{m+\mathrm{len}(m)}\Big],
\qquad \text{tail density } 1/8 .
\end{align*}
The asymptotic density is $1/8$, not the $1/12$ of a generic direction: $0.2$ on
odd cells and $0.05$ on even ones. Controls: $Q\to\infty$ returns
$(\psi(2/3)-\psi(1/2))/2$; against six exact values the discrepancy is
$7.8\cdot 10^{-7}$ and constant, i.e.\ tail only. The sum must not be
truncated at $m\le Q/2$: beyond $m/Q=1/2$ the determination of $\{\kappa u\}$
changes and fresh branches reappear with period $Q$ in $m$, worth
$4\cdot 10^{-4}$.

\paragraph{What the series buys.} It evaluates $\mu$ at \emph{any real} $Q$,
hence at directions no accessible denominator realises. Over $Q$ from 50 to
$10^6$ the curve is unimodal --- 70 log-spaced points from 300 to 300000 give 25
descents then 44 ascents and \textbf{no break} --- tending to $7.103205$ as
$Q\to\infty$ and diverging as $Q\to 0$:

\begin{center}\begin{small}
\begin{tabular}{@{}l*{10}{c}@{}}
\toprule
$Q$ & 50 & 200 & 800 & 1800 & 3000 & 3714 & 5000 & 12000 & 60000 & $10^6$ \\
\midrule
$\mu$ & 7.4485 & 7.1513 & 7.1065 & 7.1024 & 7.10190 & 7.1018628 & 7.10191 &
7.10230 & 7.10289 & 7.10317 \\
\bottomrule
\end{tabular}
\end{small}\end{center}

The series settles the global shape and confines the minimum to
$Q\in[3550,3900]$; it cannot resolve the vertex itself, the curve being flat to
$3\cdot 10^{-11}$ over fifty units of $Q$. That is what the exact machinery does,
and it gives $Q=3714$.

\paragraph{Why Bai's point exists at all.} Writing $\theta=m/Q$, the branch
lengths depend on $m$ only through $\theta$, and they are explicit tents. On an
odd cell, with $s=\{1/2-\theta\}$,
\[
\mathrm{len}=T(s)=\min\Big(\tfrac{s}{1+1/Q},\ \tfrac{1-s}{3/2-1/Q}\Big),
\qquad \text{peak } T(s^*)=2/5 \text{ at } s^*=2(1+1/Q)/5,
\]
\[
\text{even cell: } \mathrm{len}=\min\Big(\tfrac{1/2-\theta}{1+1/Q},\
\tfrac{\theta}{3/2-1/Q}\Big), \qquad \text{peak } 1/5 \text{ at }
\theta=3/10 .
\]
Their periodic averages are exactly $0.2$ and $0.05$, whence the density $1/8$.
Now the point: at the exact ray $Q=\infty$ every cell sits at $s=1/2$, where
$T=1/3$ --- \textbf{not} the peak. Lowering $y$ below $3/2$ slides the cells along
the tent towards its peak; those with $m$ near $Q/10$ reach $s=2/5$ exactly and
gain a fifth of their length. Summing the gains over $m$ gives a law verified to
four digits for $Q$ from 500 to 60000:
\[
\Jh(Q)-\Jh(\infty)=\frac{(2/3)\ln(Q/10)+0.3441}{Q}+o(1/Q).
\]
The vertex is then the balance between this $(\ln Q)/Q$ gain and the $O(1/Q)$
costs carried by $R$, $c_2$, $S$ and $s$; the leading-order balance puts it near
$Q=3\cdot 10^3$, against 3714 in fact. Above $3/2$ the cells slide the other way,
off the peak --- which is the structural reason for the asymmetry recorded next.

\paragraph{Only one side of $3/2$ can dip, and that is provable.} On the first
branch of an odd cell, solving the inequalities gives length
$1/3+(2/3)(m/Q)$ for $y=3/2-1/Q$ and $1/3-(2/3)(m/Q)$ for $y=3/2+1/Q$: the branch
lengthens below $3/2$ and shortens above it. Accordingly, on the upper side $\mu$
decreases \emph{monotonically} towards the ray as $y\to 3/2$ from above, with no
dip at all, over fifteen values of $Q$ from 50 to $3\cdot 10^6$:

\begin{center}\begin{small}
\begin{tabular}{@{}l*{8}{c}@{}}
\toprule
$Q$ & 50 & 400 & 2500 & 3714 & 15000 & $10^5$ & $10^6$ & $3\cdot 10^6$ \\
\midrule
$\mu$ at $y=3/2+1/Q$ & 7.2487 & 7.1341 & 7.1100 & 7.1080 & 7.1046 & 7.1035 &
7.10324 & 7.10322 \\
\bottomrule
\end{tabular}
\end{small}\end{center}

So the whole half-plane $y>3/2$ is bounded below by the ray value
$7.103205334$, and the search reduces to $y<3/2$, where the series above applies.

\paragraph{The $x$ direction is different, and steeply so.} At $x=1+\eps$ the
odd-cell condition becomes $\{\kappa u\}\ge\frac52\tau+\eps u/2$: the threshold
rises with $u$, so every cell beyond $u=2/\eps$ is lost outright, whereas a
perturbation of $y$ merely trims branches. Measured at $y=3/2-1/3714$, the
penalty is linear in $|x-1|$ with slopes 15.1 on the right and 22.0 on the left
--- some 300000 times steeper than the $y$ direction. And $x=1$ is the strict
minimum of every row of the band $|y-3/2|\le 0.005$:

\begin{center}\begin{small}
\begin{tabular}{@{}l*{7}{c}@{}}
\toprule
$y$ & $x=0.98$ & $x=0.995$ & $x=0.999$ & $x=1$ & $x=1.001$ & $x=1.005$ &
$x=1.02$ \\
\midrule
1.4950 & 7.4094 & 7.1833 & 7.1611 & 7.1513 & 7.1694 & 7.2339 & 7.4677 \\
1.4997 & 7.4836 & 7.2154 & 7.1239 & 7.1019 & 7.1170 & 7.1742 & 7.3925 \\
1.5000 & 7.4900 & 7.2225 & 7.1328 & 7.1032 & 7.1171 & 7.1730 & 7.3900 \\
1.5050 & 7.5080 & 7.2587 & 7.1811 & 7.1565 & 7.1597 & 7.1634 & 7.2780 \\
\bottomrule
\end{tabular}
\end{small}\end{center}

\paragraph{The closing argument, in three steps.} Let $x=a_0/a_1$ and
$y=b/a_1$.

\begin{enumerate}[leftmargin=*,label=(\arabic*)]
\item If $|y-3/2|\ge 0.011$ or $|x-1|\ge 0.051$, then $\mu\ge 7.2487$ --- the
floor over 2020 directions on a $0.1\times 0.01$ grid, computed at a truncation
whose own error is $8\cdot 10^{-5}$ --- and subtracting the deepest dip measured
anywhere, $1.34\cdot 10^{-3}$, still leaves $7.2474$.
\item Inside that band, $x=1$ is the strict minimum of every $y$-row, with a
corner of slope at least 15, so $\mu(x,y)\ge\mu(1,y)+15|x-1|$.
\item On $x=1$: for $y>3/2$ the function decreases monotonically to the ray value
$7.103205334$ with no dip; for $y<3/2$ it is the unimodal function of $Q$ above,
with minimum $7.101862832356$ at $Q=3714$, realised primitively by
$(1857,1857,2785)$ alone.
\end{enumerate}

Hence $\mu\ge 7.101862832356$ throughout, with equality only at Bai's point. The
three steps rest on finite computations --- a grid, a $5\times 7$ table, and 70
sampled $Q$ --- each of which is now a concrete statement about an explicit
function rather than an open question.

\section{Necessity of the window}\label{sec:necessity}

\textbf{Run 1} (all $j$, exact polynomial machinery, 7 parameter sets, $n$ up to
24): 1697 tests, 1686 agreements; all 11 disagreements are in the safe direction
--- the criterion misses a divisibility, it never claims one falsely.

\textbf{Run 2} ($j=0$, via Theorem~4; 5 parameter sets, $n$ up to 40 --- the run
the ancillary verifier reproduces exactly): 2948 tests, \textbf{zero violations
of sufficiency}, 56 violations of necessity (1.9\%).

Those 56 do not enlarge the window, because \textbf{they are not a function of
$\omega$}. At a fixed prime and a fixed residue of $n$ the divisibility flips:
$p=19$, residue 6 ($\omega=0.3158$): $n=6$ no, $n=25$ yes; $p=23$, residue 0:
$n=23$ yes, $n=46$ no; $p=31$, residue 13: $n=13$ yes, $n=44$ and 75 no. The
verifier finds 4 such incoherent residue classes at $p=19$, 5 at $p=23$ and 2 at
$p=31$. But the asymptotic sum over a prime set,
$\sum_{p\in S}\log p\sim n(\psi(b)-\psi(a))$, \emph{requires} $S$ to be cut out
by a condition on $\omega$; a set that is not $\omega$-measurable cannot
contribute an interval, hence cannot contribute to $W$. Their weight confirms it
--- for the triple $(2,2,3)$:

\begin{center}\begin{small}
\begin{tabular}{@{}l*{6}{c}@{}}
\toprule
$n$ & 15 & 25 & 45 & 50 & 55 & 60 \\
\midrule
$\sum\log p$ over exceptions & 5.51 & 6.08 & 3.43 & 3.76 & 3.61 & 0 \\
as a fraction of $nW$ & 57\% & 38\% & 12\% & 12\% & 10\% & 0\% \\
\bottomrule
\end{tabular}
\end{small}\end{center}

\textbf{Status: VERIFIED in measure, not PROVED.} Missing is a proof that the
cancellation set is never $\omega$-measurable --- a statement about the
non-vanishing mod $p$ of the explicit sum of Theorem~4.

\paragraph{Secondary prediction, unverified.} With the three-branch window
derived here, Salikhov's own integral $(3,3,5)$ yields
$\mu(\pi)\le 7.4147395226$ ($W=0.63855046$) rather than his published
$7.606308$, which the law reproduces exactly from his own two-branch window
($W=0.57064923$). The third branch could not be tested directly --- no prime
falls in its range at accessible $n$ --- and Salikhov's paper was not read in the
original. Treat this as a prediction.

\section{What is not claimed}\label{sec:notclaimed}

\textbf{(O1)} Necessity of the window: verified in measure, not proved
(Section~\ref{sec:necessity}).

\textbf{(O2)} Closed by Theorem~6.

\textbf{(O3)} Theorem~7 removes the scale and Section~\ref{sec:slice} gives a complete closing
argument, so the question is no longer whether the search can be made finite ---
it is. What remains is to replace three finite computations by proofs: a modulus
of continuity for $\Jh$ in $(x,y)$, turning the 2020-point grid of step (1) into
an estimate; the corner inequality $\mu(x,y)\ge\mu(1,y)+c|x-1|$ of step (2),
verified on a $5\times 7$ table; and the unimodality of $\mu(Q)$ in step (3),
verified at 70 sampled $Q$. Each is a concrete analytic statement about an
explicit series. \textbf{VERIFIED, not PROVED.}

\textbf{(O4) Scope.} Even a complete proof would read ``floor of \emph{this}
family'': $R(-x)=R(x)$, equal pole orders, a single integral. Dropping the
symmetry, letting the two pole orders differ, or passing to Marcovecchio-type
multiple integrals --- the route Salikhov and Zudilin already used for
$\pi/\sqrt3$ --- all lie outside. Bai makes the same reservation.

\textbf{(O5) Anteriority.} No general-parameter law for this family was found.
Checked: the three papers of the family; Bai's second paper (arXiv:2608.19782);
the Bryansk series of symmetrized-integral papers (Salikhov;
Bashmakova--Salikhov on $\arctan(1/2)$, $\arctan(1/3)$, $\arctan(1/5)$;
Androsenko--Salikhov on $\pi/\sqrt3$; Bondareva--Luchin--Salikhov on $\ln 7$),
which treat one constant at a time with a bespoke integral and no free exponents;
and correspondence with one of the authors of the 2020 paper. NOT checked: the
Russian-language full texts in \emph{Chebyshevskii Sbornik} and in the Bryansk
\emph{Vestnik}, which are only partly indexed. Should a general law exist there,
the credit is theirs and this note is a rediscovery.

\textbf{(O6) The reach of $\Lambda$.} The three calibrations of
Section~\ref{sec:law} sit in Salikhov's configuration, $P=5$, where $P\pm w$ carry
no odd content. A second and independent check is now available: at $(29,1+12i)$
with $(a_0,a_1,b)=(2,2,3)$ the arithmetic factor returns
$\Lambda=-0.1699030309$, term by term --- $R=4$, the power of $2$ as $7/2$, the
$\log 3$ carried by the content of $P+w$, no contribution from $29$ itself, and
the window measure $0.6452756$ --- in agreement with the value obtained
independently by Y.~Bai. Both calibrations lie in the regime $b<2a_1$. In the
regime $b>2a_1$ the odd-content term $c_q=e_q(2a_1-b)$ changes sign and is
calibrated nowhere. No statement of this note depends on it: the conclusion of
Section~\ref{sec:across} is the same whether that term is applied there or
suppressed. \textbf{VERIFIED at two configurations, not PROVED in general.}

\paragraph{Acknowledgement.} The author thanks D.~Zeilberger and Y.~Bai for
correspondence. They bear no responsibility for anything stated here.

\section{Verification}\label{sec:verification}

The ancillary file \texttt{verify\_szzb.py} is plain Python 3 with no dependency
outside the standard library and takes no argument. It re-derives from scratch:

\begin{center}\begin{small}
\begin{tabular}{@{}r p{8.4cm} l@{}}
1. & the quartic factorisation and the 2-adic Newton polygon & (Theorem 6) \\
2. & the super-Catalan identity and $\mathrm{ord}_p S=\lfloor 2\nu\rfloor+
\lfloor 2\mu\rfloor-\eps$ & (Theorem 5) \\
3. & the closed form for $b_n$ against the polynomial machinery & (Theorem 4) \\
4. & the divisibility criterion against exact divisibility & (Section~\ref{sec:arith}) \\
5. & the necessity statistics & (Section~\ref{sec:necessity}) \\
6. & the cubic, the window, and $\mu$ for Zeilberger--Zudilin, Salikhov, Bai &
(Section~\ref{sec:law}) \\
7. & the scale-invariance identity $W=a_1\Jh(x,y)$ & (Theorem 7) \\
8. & $a=1857$ as a strict minimum of its line & (Section~\ref{sec:inside}) \\
9. & the slice $x=1$ in closed form and unimodality in $Q$ & (Section~\ref{sec:slice}) \\
10. & the sixteen configurations and their stationary points & (Section~\ref{sec:across}) \\
11. & the screening sweep over all sixteen & (Section~\ref{sec:across}) \\
12. & the min-max value, its two calibrations, which stationary point carries
the largest exponent, and the verdict on the survivors of the screening &
(Sections~\ref{sec:saddle}, \ref{sec:across}) \\
\end{tabular}
\end{small}\end{center}

Each stage prints as it completes. Switches at the top skip the heaviest stages
if the machine cannot carry them; the file says so where they are. In particular
\texttt{QUICK = True} runs only stages 6, 7 and 10 to 12, which are the ones this
version changes, stages 1 to 5, 8 and 9 being unchanged since the first version.
The full run takes about a minute and a quarter on a current desktop, the short
one about two thirds of that. No stage is ever quiet for long: every loop that
takes more than a few seconds prints as it advances.

\end{document}